\documentclass[12pt]{article}
\usepackage{geometry}                
\usepackage{graphicx}
\usepackage{amssymb}					
\usepackage{amsmath}					
\usepackage{mathrsfs}
\usepackage{tikz}
\usepackage{epstopdf}
\usepackage{amsthm}
\newtheorem*{theorem*}{Theorem}
\newtheorem*{lemma*}{Lemma}
\newtheorem*{example*}{Example}
\newtheorem*{question*}{Question}
\newtheorem*{conjecture*}{Conjecture}
\newtheorem{theorem}{Theorem}
\newtheorem{lemma}[theorem]{Lemma}
\newtheorem{corollary}[theorem]{Corollary}

\newtheorem{conjecture}{Conjecture}

\def\b0{{\bf 0}}
\def\b1{{\bf 1}}

\def\cS{{\cal S}}

\def\mbr{{\mathbb{R}}}

\def\n{\noindent}

\begin{document}
\title{
Skewness, crossing number and Euler's bound for graphs on surfaces
\thanks{Submitted for publication, January 4, 2025}
}
\author{
Paul C. Kainen\\
 \texttt{kainen@georgetown.edu}}
\date{}
\newcommand{\Addresses}{{
  \bigskip
  \footnotesize

\n

\n

\par\nopagebreak
}}

\maketitle

\abstract{

\n
For every connected graph $G$ and surface $S$, we consider the well-known string of inequalities $\delta_S(G) \leq \mu_S(G)  \leq \nu_S(G)$, where $\mu$ and $\nu$ denote skewness and crossing number and $\delta$ is the Euler-formula lower bound.  Recent developments are surveyed; new results are given for the ``folded'' cube including its genus.\\
\n
{\bf Key Phrases}: 
genus, skewness, Euler lower bound, stable crossing number}.

\section{Introduction}

Several invariants measure the topological complexity of a graph;  
genus, crossing number, and skewness, in particular, are widely studied. Combinatorial invariants, based on the Euler polyhedral formula, provide lower bounds.  

This paper extends recent progress by Sun \cite{sun} on two old conjectures, due to Guy \cite{guy} and Kainen \cite{pck-lowerbound}, both published in 1972, that relate the above invariants.  A number of other recent results are unified in this framework.

The {\it genus} $\gamma$ of a graph $G=(V,E)$ is the least number of handles needed in an orientable surface in order to be able to embed $G$.
The {\it crossing number} $\nu := \nu_S$ of $G$ in a surface $S$ is the least number of crossings with which $G$ can be drawn in the surface, while
{\it skewness} $\mu:= \mu_S$ is the least number of edges such that deleting the edges makes $G$ embeddable in $S$.  Clearly, $\mu$ is a lower bound on crossing number, $\mu \leq \nu$.  

When the surface has at least as many handles as the genus of $G$, then crossing number and skewness are zero, and conversely.  As with {\it cycle rank}, genus, skewness, and crossing number (Harary \cite[pp 39, 116--119, 122--124]{harary}) are not changed by edge subdivision; that is, they are {\it topological} invariants.

The Euler {\it excess} $\delta := \delta_t$ of connected $G$ in orientable surface $S_t$ \cite{pck-lowerbound} is 
\begin{equation}
\delta_t(G) = q - \alpha(p - 2 + 2t),
\label{eq:delta}
\end{equation}
where $q = |E_G|$ and $p = |V_G|$, $t$ is the number of handles in the surface, and $\alpha = \frac{g}{g-2}$ with $g$ equal to {\it girth}  (length of shortest cycle).  
Excess 
is {\it metric}.  

For $G$ acyclic, we set $\delta(G) = 0$ as forests are always embeddable.

One has $\delta \leq \mu$ since the subtracted expression $\alpha(p - 2 + 2t)$ in (\ref{eq:delta}) is the largest possible number of edges in any $p$-vertex $S_t$-graph with girth $g$, and further $\delta = \mu$ exactly when it is possible to remove a subset $E' \subset E(G)$ of edges such that $G - E'$ embeds in the surface and has maximum edge-count.

Set
$\varepsilon := \max\{0, \lceil \delta \rceil \}$;
we call $\varepsilon$ the {\it integer excess}.  As $\mu$ is non-negative and integer-valued, it follows that 
\begin{equation}
\delta \leq \varepsilon \leq \mu \leq \nu,
\end{equation}
and {\it it is natural to wonder about equalities and differences}.

Guy \cite{guy} showed that for fixed $t$, $\delta_t(G) = \mu_t(G)$ for $G$ of sufficiently large order if $G = K_{a,b}$, $a, b$ even, and $G = Q_d$ and he conjectured it for all complete and bipartite complete graphs.  In contrast, we conjectured \cite{pck-lowerbound} $\delta_t(G) = \nu_t(G)$ for $G$ complete or bipartite complete 
when $t$ is just smaller than the genus (a stronger claim for a smaller domain), and we proved \cite{pck-stable} for $0 \leq k \leq 2^{d-4}$,
 $4k = \delta_t(Q_d)= \mu_t(Q_d) \leq \nu_t(Q_d) \leq 8k$ when $t = \gamma(Q_d) - k$.  

In the following sections, we review recent progress on these problems.
Note that in general, {\it crossing number is not near its Euler lower bound}.  In the plane, for the bipartite graph $K_{a,b}$, the lower bound is $$\delta_0(K_{a,b}) = ab - 2(a+b-2) = (a-2)(b-2),$$ while $\nu_0(K_{a,b}) \sim c a^2 b^2/16$ for some $c$, $.91 < c \leq 1$ (Balogh et al \cite{balogh}).  But for all $a, b$, {\it planar skewness and excess are equal for $K_{a,b}$} as there are quadrangulations of the plane that span $K_{a,b}$.

If $G$ has a cycle, then $\delta$ becomes negative once $t$ is large enough. Define the {\it algebraic genus} of $G$, $h(G)$, to be the maximum $t$ such that $\delta_t(G) \geq 0$.
In \cite{pck-lowerbound} we conjectured that the equation \framebox{$\delta = \nu$} holds for the pairs $(K_n, S_t)$ when $t = h(G)$.
 Riskin \cite{riskin} showed that, for $K_9$ in $S_2$, this claim is off by 1, but Sun \cite{sun, sun-abs} showed that otherwise the conjecture is correct and, suitably modified, that it also holds for non-orientable surfaces except for well-known examples $K_7$ due to Franklin \cite{franklin} and $K_8$, due to Ringel \cite{ringel}.

The corresponding $\delta = \nu$ conjecture from \cite{pck-lowerbound}, for  complete {\it bipartite} graphs, is proved in \cite{sun} with $K_{3,5}$ and $K_{5,5}$ excepted, where the difference is 1.  But {\it skewness equals its lower bound for the two exceptional graphs} as we note below.  

Guy \cite{guy} gave an elegant proof that \framebox{$\delta = \mu$} for some bipartite complete graphs.  In Section 4 below, we describe a possible extension of his method for $K_n$. This genus technique might be added to those described in \cite{wan}.


One may prefer skewness to crossing number for various reasons.
Even in the plane, crossing numbers are notoriously difficult to calculate, while $t$-skewness is more tractable.  This has led to a growing literature on properties of skewness (e.g., \cite{bokal, chia, cim, cook, farr, ng, pck-stable, pck-5ct, pck-skewness, liebers, liu-pc, ouyang, ouyang2, pahlings, planariz}).

In applications, removing edges from the graph is often the natural way to avoid crossings.  For instance, to avoid contact
between two electrical conducting paths in a VLSI circuit, one must fabricate an ``overpass'' (called a {\it via}), increasing cost and decreasing reliability.

Deleting a single edge can delete multiple crossings; in fact, starting with a 3-connected maximal planar graph, one can add a single edge to obtain a graph with planar skewness 1 but arbitrarily large crossing number.

We show that the fixed graphs $K_n$ and $K_{a,b}$ satisfy the equation \framebox{$\varepsilon = \mu$} for the surface $S_t$ when $t$ is small enough or large enough and for an intermediate value of $t$. 
Further, we show that for both hypercubes and ``folded cubes'' \cite{el-A-L}, the equation $\varepsilon = \mu$ holds for {\it every} orientable surface, and the genus of folded cubes is determined, cf \cite{hao-cair}.  



An outline follows.  Section 2 has definitions and lemmas; planar and toroidal cases for complete and bipartite complete graphs are in Section 3, while Section 4 on orientable surfaces of higher genus extends Sun's work \cite{sun} and Guy's \cite{guy}.  Cubes and folded cubes, in Section 5, are shown to satisfy the equation $\varepsilon = \mu$ for {\it all} surfaces.
Section 6 is a discussion.

\section{Definitions and lemmas}

Graphs $G=(V,E)$ are simple; terminology is from \cite{harary}.  Complete graphs and bipartite complete graphs are denoted $K_n$ and $K_{a,b}$, resp.  See, e.g., \cite{pck-gw, pck-amsuh}.

Put the vertices in $V$ into $\mbr^3$ as distinct points and add each edge in $E$ as the straight-line segment joining its  endpoints. It is always possible to do this such that no three vertices lie on a straight line and no two edges intersect except at a common endpoint.  This {\bf 3-d  realization} $|G|$ of $G$ inherits a topology as a subset of Euclidean 3-space and this topology does not depend on the placement; it is the unique topology that makes all edges homeomorphic to the unit interval. Non-endpoints of an edge are called {\bf interior} points.

A {\bf handle} is a topological space homeomorphic to the annulus. It is  attached to the sphere (or to a surface or disjoint union of surfaces) by deleting the interior of two disjoint closed disks contained in the sphere or surface or disjoint union of surfaces and then identifying the two boundary circles of the handle with the boundaries of the open disks that have been removed.  This operation preserves the locally Euclidean property. 
A {\bf crosscap} is a space homeomorphic to the M\"obius strip and its attachment requires the removal of a single open disk, again with identification of the boundary circles.

A {\bf surface} $S$ is a 2-manifold, a connected, Hausdorff topological space which is locally Euclidean (each point has a neighborhood homeomorphic to $\mbr^2$).  
There are two kinds.  The {\bf orientable} surfaces $S_t$ consist of a sphere with $t \geq 0$ handles attached; the {\bf non-orientable} $N_c$ are spheres with $c > 0$ crosscaps. We focus on the orientable case.

An {\bf embedding} of graph $G$ in surface $S$ is a continuous 1-to-1 function $f:|G| \to S$.  We usually write $f:G \to S$.  This function is a homeomorphism onto its image $f(G)$ which is compact and hence closed.  The {\bf regions} of the embedding are the open path-connected components of the complement $S \setminus f(G)$ of $f(G)$.

An embedding is {\bf 2-cell} ({\bf closed 2-cell}) if each region is an
open disk (or if the {\it closure} of each region is a closed disk). Let $f: G < S$ and $f: G \leq S$ denote that $f$ is an open (or closed) 2-cell embedding of $G$ in surface $S$.
The {\bf Euler characteristic} $\chi(S)$ is $2-2t$ if $S=S_t$, while $\chi(N_c) = 2 - c$.  

For a 2-cell embedding of connected graph $G$ in $S$, {\bf Euler's formula} says
\begin{equation}
p - q + r = \chi(S),
\end{equation}
where $p$ and $q$ count vertices and edges and $r$ counts regions in the embedding; $p$ is the {\bf order} and $q$ the {\bf size} of the graph.

A {\bf drawing} of a graph $G=(V,E)$ in a surface $S$ is a continuous map $f: |G| \to S$ which is 1-to-1 restricted to the set of vertices and is also 1-to-1 restricted to each separate edge $e \in E$ in the 3-d realization.  Thus, the nodes of the graph are points in the surface and each edge $e = vw$  of the graph is represented in the surface by a non-singular curve $f(e)$ joining $f(v)$ and $f(w)$.  

A {\bf crossing} in a drawing $f$ is a distinct pair of interior points $x, x'$ in distinct edges such that $f(x) = f(x')$.  
We shall {\it require} that drawings satisfy: {\it if $y$ is in two distinct edges, then $y$ is either a crossing or a common endpoint}. Thus, a drawing is an embedding if and only if it is crossing-free.
 
The {\bf genus} $\gamma(G)$ of a graph $G$ is the least $t$ such that $G$ can be embedded in $S_t$.  An embedding of $G$ in $S_t$ with minimum $t$ is called a {\bf minimum-genus embedding}; such an embedding is necessarily 2-cell \cite{youngs}.  A 2-cell embedding of connected $G$ is minimum genus if and only if it has the greatest possible number of regions over all 2-cell embeddings of $G$ in any orientable surface. 

An embedding of $G$ in $S$ is a {\bf triangulation} if each region has 3 sides; an embedding is a {\bf quadrangulation} (a {\bf quad embedding}) if each region has 4 sides.  The graph $G$ is called the {\bf skeleton} of the tri- or quadrangulation. An $n$-vertex triangulation of $S_t$ has skeleton a maximum-size subgraph of $K_n$, while an $n$-vertex bipartite quadrangulation has skeleton a maximum-size subgraph of $K_{a,b}$ if $G$ is bipartite with a bipartition having cardinalities of $a$ and $b$.  

For $t \geq 0$, the $t$-{\bf crossing number} $\nu_t(G)$ is  the least number of crossings in a drawing of $G = (V,E)$ in $S_t$.
The $t$-{\bf skewness} $\mu_t(G)$ (Guy \cite{guy} called it ``the generic slimming number'', see also \cite{pck-lowerbound} and Pahlings \cite{pahlings}) is the minimum of $|E'|$, where $E' \subset E$ and $G - E'$ embeds in $S_t$.   Harary \cite[p 124]{harary}, 
cites Kotzig for the original concept of planar skewness, see \cite{kotzig55, rosa}.

If $G$ is a forest, define $\delta_t(G) := 0$ for all $t$, while if $G$ is connected but not a tree,
 for $t \geq 0$, put
\[
\delta_t(G) := |E_G| - \Big(\frac{g}{g-2}\Big)(p - 2 + 2t),
\]
where $g := g(G)$ is the girth; we call $\delta_t$ the $t$-{\bf excess}. The {\bf integer $t$-excess} is
\begin{equation}
\varepsilon_t(G) := \max \{0, \lceil \delta_t(G) \rceil \}.
\end{equation}



\begin{lemma}[\cite{pck-lowerbound}]
For all graphs $G$ and $t \geq 0$, 
$\delta_t(G) \leq \varepsilon_t(G) \leq \mu_t(G) \leq \nu_t(G)$.
\end{lemma}

For any graph $G$, 
skewness and excess are equal
for all but a finite set of orientable surfaces.

\begin{lemma}
For all graphs $G$, if $t \geq \gamma(G)$, then $\varepsilon_t(G) = 0 = \nu_t(G) = \mu_t(G)$.
\label{lm:all-but}
\end{lemma}
We first deal with complete graphs and bipartite complete graphs, where {\it skewness equals excess with specific exceptions}; then we show that the equality {\it always} holds for the cube and its supergraph, the folded cube.


\section{$K_n$ and $K_{a,b}$ in the torus and the sphere}

We consider the sphere and torus, the two orientable surfaces of lowest genus, and we show that skewness minus excess is small for two fundamental families.


\begin{theorem}
If $n \geq 1$, then $\varepsilon_0(K_n) = \mu_0(K_n)$ 
\end{theorem}
\begin{proof}
If $n \leq 4$, $K_n$ has genus zero so  Lemma \ref{lm:all-but} suffices. If $n \geq 3$, there is an $n$-vertex triangulation of $S_0$, which embeds a spanning subgraph of $K_n$.
\end{proof}

\begin{theorem}
If $a \geq b \geq 1$, then $\varepsilon_0(K_{a,b}) = \mu_0(K_{a,b})$.
\label{th:planar-bip}
\end{theorem}
\begin{proof}
If $b = \min\{a,b\} = 1$, $K_{a,b}$ is planar so Lemma \ref{lm:all-but} applies.  If $b \geq 2$, there is a bipartite quadrangulation of $S_0$ with $a$ vertices in one part and $b$ in the other part. Indeed, for $a = 2 = b$, $C_4 = K_{2,2}$ works.  Given any bipartite quadrangulation of $S_0$ with bipartition having parts of size $a'$ and $b'$, one can add a new vertex in the interior of one of the $C_4$'s and join it to the two vertices in the $a'$ part to get a larger bipartition quadrangulation of $S_0$. 
\end{proof}

Similarly, excess equals skewness for complete graphs on the torus.

\begin{theorem}
Let $n \geq 1$.  Then $\varepsilon_1(K_n) = \mu_1(K_n)$.
\end{theorem}
\begin{proof}
If $n \leq 7$, then $K_n$ has genus not exceeding $1$ so  use Lemma \ref{lm:all-but}. 
 
If $n = 2r$, $r \geq 4$, embed $r$-cycles $u_1, \ldots, u_r$  and $w_1, \ldots, w_r$ horizontally along the top and bottom of the standard torus in 3-space, where the torus is oriented so that the xy-plane intersects it in two concentric circles, while both the xz and yz planes intersect the torus in two disjoint, nonconcentric circles.  
For each $j \in [n]$, one can put two edges on the front face, making the vertex $u_j$ in position $j$ on the top $u$-cycle adjacent to the nodes $w_j$ and $w_{j+1}$ on the bottom cycle (addition mod $r$).  Also, one can put two edges on the back face, joining position $w_j$ to positions $u_{j+1}$ and $u_{j+2}$ on the top.  Each vertex has 
degree 6 and each region is a triangle.  Hence, the embedding is minimum.

This suffices but we verify the count explicitly. In $K_{2r}$, there are $r^2$ edges between the two $r$ cycles and we only have used $4r$ of them so $r (r-4)$ are deleted. For the edges between vertices in the same $r$ cycle, we only 
include $r$, corresponding to the consecutive pairs, of the total number ${r \choose 2} = |E(K_r)|$ so the total number of deleted edges is $r^2 - 4r + r(r-1) - 2r = 2r^2 - 7r$.
Hence, $\mu_1(K_{2r}) \leq 2r^2 - 7r$.  But $\mu_1(K_{2r}) \geq \delta_1(K_{2r}) = 2r^2 - 7r$.

If $n = 2r+1$, $r \geq 4$, first put $K_{2r}$ in $S_1$ by deleting $2r^2 - 7r$ edges.  Then put the $2r+1$-st vertex into the interior of one of the triangles and join it to the three vertices of the triangle, deleting the $2r-3$ additional edges.  Hence, $\mu_1(K_{2r+1}) \leq 2r^2 - 7r + 2r -3 = 2r^2 - 5r -3 = \delta_1(K_{2r+1}) \leq \mu_1(K_{2r+1})$.  
\end{proof}

In the bipartite complete case, for the torus, the difference is at most 2.
\begin{theorem}
Let $a \geq b \geq 1$. Then $\varepsilon_1(K_{a,b}) = \mu_1(K_{a,b})$ unless $a \geq 7$ and $b=3$, in which case $\delta_1(K_{a,3}) = a-6 \leq \mu_1(K_{a,3}) \leq a-4$.  
\label{th:quad-nest}
\end{theorem}
\begin{proof}
If 
$b \leq a \leq 6$, then $K_{a,b}$ is toroidal so use Lemma \ref{lm:all-but}.  If $\min \{a,b\} \geq 4$, quadrangulate $S_1$ by $K_{4,4}$ and proceed as in the proof of Theorem \ref{th:planar-bip}.  For $a \geq 7$, an upper bound on $\mu_1(K_{a,3})$ is obtained by extending the planar drawing of $K_{a,2}$ with $a$ points on the $y$-axis and two points at the points $\pm 1$ on the $x$-axis.  The third point is now placed at $+2$ on the $x$-axis and joined to the two outermost $y$-axis points; a handle is added joining innermost and outermost quadrilaterals in the planar drawing of $K_{a,2}$, and the third point is joined to the two innermost $y$-axis points using the handle.  Thus, the added point has degree 4 so $a-4$ of its incident edges are deleted.
\end{proof}

\section{$K_n$ and $K_{a,b}$ on higher genus surfaces}

Given a graph $G$, we ask: 
{\it For which $t \geq 0$ is $\delta_{t}(G) = \mu_t {G}$}?  
By the previous section, we know that for $n \geq 1$, equality holds for $G = K_n$ if $t \leq 1$, or, by the Ringel-Youngs genus formula for $K_n$ \cite[p 120]{harary} and Lemma \ref{lm:all-but}, if
$$t \geq \gamma(K_n) :=\Big\lceil \frac{(n-3)(n-4)}{12}  \Big\rceil.$$  
Also for $a \geq b \geq 1$ if $t \leq 1$ or, by Ringel's genus formula for $K_{a,b}$  \cite{harary}, if 
$$t \geq \gamma(K_{a,b}): = \Big \lceil \frac{(a-2)(b-2)}{4} \Big \rceil.$$ 

Thus, for $G$ complete or bipartite complete, the open cases for orientable surfaces $S_t$ is $2 \leq t \leq \gamma(G)-1$.  We now consider these.

The {\bf algebraic genus} $h := h(G)$
of a graph $G$ is $h = \max\{t: \delta_t(G) \geq 0\}$; so $h(K_n) = \lfloor \frac{(n-3)(n-4)}{12}  \rfloor$ and $h(K_{a,b}) = \lfloor \frac{(a-2(b-2)}{4}  \rfloor$.
We conjectured in \cite{pck-lowerbound} that for complete or bipartite complete graphs $G$,
$\delta_h(G) = \nu_h(G)$.  Sun \cite{sun} showed the conjecture is correct except
for several cases where it is off by $\leq 1$.

\begin{theorem}[Sun]
{\rm (i)} For $n \neq 9$, $\delta_h(K_n) =  \nu_h(K_n)$ and $\nu_h(K_9) - \delta_h(K_9) = 1$.\\
{\rm (ii)} For $(a,b) \notin \{(5,3), (5,5)\}$, $\delta_h(K_{a,b}) =  \nu_h(K_{a,b})$, and $\nu_h(K_{5,b}) {-} \delta_h(K_{5,b}) {=} 1$
for $b =3$ or 5. {\rm (iii)} Equality holds in non-orientable cases except $K_7$ and $K_8$.
\label{th:Sun}
\end{theorem}

As $\delta_h \leq \mu_h \leq \nu_h$, the above cases of equality for excess and crossing number also hold for excess and skewness.  For $K_9$, where $h = 2$ and $ \delta_h(K_9) = 3$, Sun observes that a result of Huneke \cite{huneke} states there is no triangulation of $S_2$ with 9 vertices, but an embedding of $K_9 - \{e_1, e_2, e_3\}$ in $S_2$ would have to be a triangulation.  So $4 \leq \mu_2(K_9) \leq \nu_2(K_9) = 4$, the latter by \cite{riskin}. Thus, $\mu - \delta = 1$ for $K_9$ in the orientable and $K_7$ and $K_8$ in the non-orientable cases.

For the exceptional {\it bipartite} cases, one has
$\delta_0(K_{5,3}) = 15 - 2(8 - 2) = 3$ and $\delta_2(K_{5,5})=1$.  While Kleitman \cite{kleitman} showed $\nu_0(K_{5,3}) = 4$ (hence, an exceptional case for Theorem \ref{th:Sun}), it is easy to see that $\mu_0(K_{5,3}) = 3$ using the 
Zarankiewicz drawing \cite{zaran}, which places the two vertex classes on the two planar axes, with each class divided as equally as possible in the positive and negative half-axes.  Further, $\mu_2(K_{5,5}) = 1$  \cite{mohar}, \cite[Fig. 34]{sun}. Thus, $\mu_h = \delta_h$ for complete bipartite graphs with {\it no exceptional cases}, where $h$ is the algebraic genus.

\begin{corollary}
For $a \geq b \geq 1$, $\delta_h(K_{a,b}) =  \mu_h(K_{a,b})$.
\end{corollary}

Sun extended his previous work to show \cite{sun-Dec} that $\nu_t(K_n) = 6k = \delta_t(K_n)$ when $t = \gamma(K_n) - k$ provided $n = 12s$, $s \geq 1$, and $0 \leq k \leq s$, and also that  $\nu_t(K_n) = 6k+3 = \delta_t(K_n)$ when $t = \gamma(K_n) - k$ provided $n = 12s + 6$, $s \geq 0$, and $1 \leq k \leq s+1$.  These imply corresponding equalities for excess with skewness.

One might expect some sort of reasonable behavior for the function
\begin{equation}
f_n(t) := \mu_t(K_n) - \varepsilon_t(K_n)
\end{equation}
so it is interesting that when $n = 2r$ for $r \not\equiv 2$ (mod 3), one has
\[ f_n\Big(\frac{(r-1)(r-3)}{3}\Big) = 0. \]

Indeed, the octahedron $O_r$ is obtained from $K_{2r}$ by deleting $r$ pairwise disjoint copies of $K_2$ and $\gamma(O_r) = \lceil \frac{(r-1)(r-3)}{3} \rceil$. Hence, for $r \equiv 1$ or 3 (mod 3), $\mu_{t_O}(K_{2r}) \leq r$ when $t_O = \frac{(r-1)(r-3)}{3} = \gamma(K_{2r}) - \lceil \frac{r}{6} \rceil + 1$.  But 
\[ \delta_{t_O}(K_{2r}) = {2r \choose 2} - 3\Big(2r-2+2\Big(\frac{(r-1)(r-3)}{3}\Big)\Big) = r,\]
and therefore $t_O$-skewness equals $t_O$-excess 
for such $K_{2r}$.
A different intermediate result could be obtained using a result of Gross \cite{gross}.

Similarly, for $t_{Q} := \gamma(Q_d) = 1 + (d-4)2^{d-3}$, we have 
\[ \epsilon_{t_Q}(K_{2^{d-1}, \;2^{d-1}}) =  \mu_{t_Q}(K_{2^{d-1}, \; 2^{d-1}})\]
as $Q_d$ has a quadrilateral embedding in its genus surface $S_{t_Q}$ and $Q_d$ is a  spanning subgraph of $K_{2^{d-1},2^{d-1}}$.  So again skewness minus excess is zero for these regular complete bipartite graphs on surfaces with small and large genus and in the intermediate case,
$t_Q = 1 + (d-4)2^{d-3} << \gamma(K_{2^{d-1}, \; 2^{d-1}}) \sim 2^{2d-4}.$  

For $t \geq 0$, let $\tau(t)$ be the least number of vertices in a triangulation of $S_t$, \cite{jung-ringel}, \cite[A123869]{oeis};
Quadrangulations have been widely studied, e.g., \cite{abu, brown-jackson, hunter-pck, liu, ringel} but relatively little seems to have been done restricting to bipartite graphs.  
Magajna, Mohar and Pisanski \cite{maga}
found $b(S) = \lceil 4 + \sqrt{16 - 8 \chi(S)} \rceil$, where $b(S)$ denotes the least number of vertices in a bipartite quadrangulation of $S$. However, we need to control the bipartition. 
Write $\xi(t)$ for the least $r$ such that $S_t$ has a bipartite quadrangulation with bipartition $(r,r)$.  When $r$ is even and $t = (r-2)^2/4$,  $\xi(t) = r$ and $b(S_t) = 2r$.  Hence, as above, one has
\begin{theorem}
If $t \geq 0$, then (i) $\varepsilon_t(K_n) = \mu_t(K_n)$ for all $n \geq \tau(t)$, and\\
(ii) $\varepsilon_t(K_{a,b}) = \mu_t(K_{a,b})$ for all $a \ge b \geq \xi(t)$.
\end{theorem}

Guy \cite[Fig. 1]{guy} gives a drawing for a subgraph of $K_{a,b}$, $a, b$  even.  Divide the $a$ and $b$ points equally and put them at the positive and negative integer points on the $x$ and $y$ axes. Join the uppermost and lowermost $y$-axis points to all $a$ points on the $x$ axis and join the innermost pair of $x$-axis points to all $b$ points on the $y$-axis. This quadrangulates the plane with $2(a+b) - 4$ lines 
and shows that $\mu_0(K_{a,b}) = \delta_0(K_{a,b})$. 

Guy then places $(a/2)-1$ ``trees'' of  branching handles into gaps between $x$-axis points, each going to similar branching ``trees'' in $(b/2)-1$ gaps between inner $y$-axis points. In modern terminology, one of Guy's ``trees''  is a pair-of-pants but allows bifurcation into any number of branches.  The $y$-axis trees use the central gap if and only if $b \equiv 0$ (mod 4) and if $a \equiv 2$ (mod 4), the $x$-axis trees avoid the central gap and are placed between non-inner vertices.  See \cite[Fig. 4]{guy}.  The previous planar surface now has genus $(a-2)(b-2)/4$.  Further, all the edges of $K_{a,b}$ that are not in the plane are contained, 4 per handle.  For instance, as he shows, $K_{6,8}$ embeds in $S_6$ by extending the quadrangulation with $24 = 2(6 + 8 - 2)$ edges with the remaining 24 edges on the 6 handles.

We believe this construction works also for $K_n$ when $n$ satisfies suitable divisibility conditions.  For example, when $n = 6m$, one could begin with $m$ disjoint spheres in $\mbr^3$, where each sphere is triangulated by $O_3 = K_6 - 3K_2$.  Each of 2 vertex-disjoint triangular faces is the location of a pair-of-pants that branches to $2(m-1)$ pairwise-vertex-disjoint faces on the remaining copies of $K_6 \leq S_0$, and each handle carries 6 edges. The following is a consequence.

\begin{conjecture}
If $n = 6m$ for $m \geq 1$, then $\delta_t(K_n) = \mu_t(K_n)$ for
$$m-1 \leq t \leq 2m^2 - 3m +1.$$
\end{conjecture}
In fact, one should prove that for $t = 4 {m \choose 2} - (m - 1) - k = 2m^2 - 3m +1 - k$, one gets $\delta_t(K_n) = 6m^2 - 3m + 6k = \mu_t(K_n)$ when $ 0 \leq k \leq 2m^2 - 3m +1 - (m-1)$.

\section{Cubes and folded cubes}

In this section, we show that, for hypercubes, {\it skewness and Euler excess are equal on all orientable surfaces}.  Arguments are based on a natural embedding of the hypercube into three-dimensional space. The folded cube is also studied.

For $d \geq 0$, the $d$-{\bf cube} (or {\it $d$-dimensional hypercube}) $Q_d$ is the graph with $2^d$ vertices, where each vertex is a length-$d$ vector of 0s and 1s.  Two vertices are adjacent if they differ in a single coordinate; $Q_d$ is $d$-regular and bipartite.

Each cube $Q_d$ is the $d$-times-iterated Cartesian product graph of the edge graph $K_2$ so $Q_0 := K_1$ and for $d \geq 1$,
\begin{equation}
Q_d = Q_{d-1} \,\Box \,K_2.
\label{eq:cube-rec}
\end{equation}
We extend (\ref{eq:cube-rec})  to a quad embedding of $Q_d$ in $S_t \subset \mbr^3$ for $t = \gamma(Q_d)$. 

More generally, let $G$ be any triangle-free connected graph.  A {\bf vertex-disjoint quadrilateral cover VDQC} of $G \leq S_t$ is a {\it set $\cS$ of 4-sided regions} such that {\it each vertex of $G$ is in the boundary of exactly one} member of $\cS$ \cite{hunter-pck, pck-stable, pisansky}; e.g., $Q_3 \subset S_0$ has VDQC consisting of any two opposite faces. 

\begin{lemma}
If $f: G \leq S_t$ has  VDQC $\cS$, then there exists $F: G \,\Box\, K_2 \leq S_{t'}$ with $t' = 2t {+} |\cS| {-} 1$ with VDQC $\cS'$ and $|\cS'| = 2 |\cS|$.
\label{lm:prod}
\end{lemma}
\begin{proof}
Take two disjoint mirror-image identical copies $f_1$ and $f_2$ of $f$ and mirror copies $\cS_1$ and  $\cS_2$ of $\cS$ in $\mathbb{R}^3$ and attach tubes (handles) joining each quad in $\cS_1$ to the corresponding quad in $\cS_2$.  Such a  handle carries the 4 edges in $G \,\Box\, K_2$ that join the 4 vertices in the first quad to the corresponding vertices in the second quad.  The 8 vertices of $G\,\Box\, K_2$ in the two disjoint embeddings $f_1$, $f_2$, that were covered by the two quads,  are now covered by two vertex-disjoint quads in the annulus $C_4 \,\Box\, K_2$ in the resulting embedding $F$.
\end{proof}
Applying Lemma \ref{lm:prod} to $Q_3 \leq S_0$ with $\cS = 2$ gives $Q_4 \leq S_1$ as $1 = 2*0 + 2 - 1$.
\begin{corollary}
Let $d \geq 3$. Then $\gamma(Q_d) = 1 + (d-4)2^{d-3}$ with a quad embedding.
\label{co:cube-genus}
\end{corollary}
\begin{proof}
By induction from Lemma \ref{lm:prod}, $\gamma(Q_d) \leq 1 + (d{-}4)2^{d-3}$ for $d \geq 3$, with the example as basis case; equality follows from the Euler lower bound $\gamma(G) \geq 1 - p/2 + q/4$ for connected graphs with girth at least 4 \cite[p 118]{harary}.
\end{proof}
This well-known result can be obtained in other ways, e.g., by factoring the 2-skeleton of the $d$-dimensional hypercube into embedding surfaces \cite{hk-monthly}.  However, the method shown
here is simple and allows for modifiications and extensions.  For instance, we showed the following in \cite{pck-stable}.

\begin{theorem}
If $d \geq 3$, $0 \leq k \leq 2^{d-4}$, then 
$$\delta_{\gamma(Q_d) - k}(Q_d) = 4k = \mu_{\gamma(Q_d) - k}(Q_d) \leq \nu_{\gamma(Q_d) - k}(Q_d) \leq 8k.$$
\label{th:stable}
\end{theorem}

One can extend the first part of Theorem \ref{th:stable}.  Recall $\varepsilon_t := \max\{0, \lceil \delta_t \rceil\}$.
\begin{theorem}[\cite{guy,pck-lowerbound} ]
If $d \geq 1$ and $t \geq 0$, then $\mu_t(Q_d) = \varepsilon_t(Q_d)$.
\label{th:mu=eps}
\end{theorem}
\begin{proof}
If $t \geq \gamma(Q_d)$, both invariants are zero.  
For $Q_d$ in $S_t$, $0 \leq t < \gamma(Q_d)$,  with $d \geq 4$, map $Q_d$ into its genus surface $S(d) := S_{\gamma(Q_d)} \subset \mbr^3$ such that we have a $d{-}2$-cube consisting of vertex-disjoint $2$-cubes on disjoint copies of $S_0$ with pairwise disjoint handles corresponding to the $(d{-}2)2^{d-3}$ edges of the $d{-}2$-cube. The cycle rank \cite[p 39]{harary} of $Q_{d-2} = (d{-}2)2^{d-3} - 2^{d-2} + 1 = \gamma(Q_d)$.  

Note that the handles don't all join the 2-cubes directly.  For instance, $Q_4$ is a square of squares.  Tubes 1 and 2 join the two copies of $Q_2$ that make the first and second $Q_3$ in $S_0$ and then faces on these two tubes are joined by an additional pair of tubes, using the fact that each product with $K_2$ has a first and second element to choose corresponding faces, creating $Q_4$ in $S_1$. 

 If $t = \gamma(Q_d) - k$, for $1 \leq k \leq \gamma(Q_d)$, then select $k$ edges in $Q_{d-2}$, $\{\overline{e}_1, \ldots, \overline{e}_k\}$, which do not disconnect the $(d{-}2)$-cube and remove the corresponding $4k$ edges $\{(u,\overline{e}_j) : 1 \leq j \leq k, \,u \in V(Q_2)\}$ of $Q_d$ to get a subgraph $Q'_d$ of $Q_d$ with a quad embedding $Q'_d \leq S_t$.  Hence, $\mu_t(Q_d) = \varepsilon_t(Q_d)$. \end{proof}
 
 If $k=1$, a single handle is removed, while if $k = \gamma(Q_d)$, then the cycle rank of the connected graph $Q_{d-2} - \{\overline{e}_1, \ldots, \overline{e}_k\}$ is zero, and this subgraph of $Q_{d-2}$ is a spanning tree.
This gives a quadrilateral embedding of a subgraph of $Q_d$ into $S_0$.

We now extend the above results to a supergraph of the cube.

Each vertex $v$ in the cube $Q_d$ has graph distance $d$ from its reversal $\overline{v}$, where 0 and 1 are interchanged; vertex $\overline{v}$ is  {\bf the antipode} to $v$ as their distance $d(v,\overline{v}) = \max_wd(v,w) = d$ and $d(v,w)=d$ only if $w = \overline{v}$. 
Let $d \geq 2$.  The {\bf folded}  $d$-{\bf cube} $F_d$, defined in  \cite{el-A-L}, is the supergraph of $Q_d$ obtained by appending the $2^{d-1}$ edges of form $v \overline{v}$. 
For instance, $F_2 = K_4$ and $F_3 = K_{4,4}$.

Note that $F_d$ is {\it bipartite} if and only if $d \geq 3$ is {\it odd}, and for all $d \geq 3$, $F_d$ has girth 4.  
Our model for the cube in 3-space extends to the folded cube. 
\begin{theorem}
Let $d \geq 3$. Then $\gamma(F_d) = 1 + (d-3)2^{d-3} = \gamma(Q_d) + 2^{d-3}$.
\end{theorem}
\begin{proof}
The lower bound is from \cite[p 118]{harary} as in the proof of Corollary \ref{co:cube-genus}.
For the upper bound we show that the antipodal edges $\{v\overline{v}: v \in V(Q_{d-1})\}$ can be carried (4 per handle) by an additional family of $2^{d-3}$ handles.  

This is clear for $F_3$ which consists of $Q_3$ on a sphere in $\mbr^3$ with a handle joining the square with last coordinate 0 to the square with last coordinate 1; the 4 non-cube edges are the antipodal edges and are carried by the handle.  

Recursively having embedded $F_{d-1}$ as $Q_{d-1}$ with  $2^{d-4}$ additional handles, make two mirror-image copies of $F_{d-1}$ in $\mbr^3$.  Each of the two copies $X_0$ and $X_1$ of $Q_{d-1} \subset S_{\gamma(Q_{d-1})}$ has a canonically defined VDQC, $\cS_0$ and $\cS_1$.  

By adding handles that join the corresponding quads in the VDQC's, one gets a genus embedding of $Q_d$, where $X_0$ and $X_1$ correspond to the points with $d$-th coordinate equal 0 or 1, resp.
Connecting the two disjoint copies of $F_{d-1}$ by the above $2^{d-3}$ handles gives a surface of the correct genus with a quad embedding of  $F_{d-1} \,\Box\, K_2$.  We perform a {\it multi-component cobordism} (a surgery involving multiple components) by switching the ends of corresponding antipodal handles to create $F_d$. More exactly,
for each pair of corresponding handles that carry the antipodal edges in the two embedded copies of $F_{d-1}$, use the mirror-image isomorphism to fix the first end of the handles but interchange their other ends.  The genus is left unchanged.
\end{proof}

To prove the analog of Theorem \ref{th:mu=eps}, we make a few observations.  For $d \geq 2$, $F_d$ has $2^{d-1}$ edges that are not in $Q_d$.  For $d \geq 3$, $F_d$ and $Q_d$ have the same girth and the same number of vertices, so for $t \geq 0$, $\delta_t(F_d) - \delta_t(Q_d) = 2^{d-1}$.
\begin{theorem}
If $d \geq 2$ and $t \geq 0$, then $\mu_t(F_d) = \varepsilon_t(F_d)$.
\end{theorem}
\begin{proof}
For $d = 2$, the result is trivial.  For $d \geq 3$, if $t \geq \gamma(F_d)$, both skewness and excess are zero.  For $t < \gamma(F_d)$, by Theorem \ref{th:mu=eps} and the above remarks, it suffices to show that $\mu_t(F_d) - \mu_t(Q_d) = 2^{d-1}$.  Clearly, the difference can't exceed $2^{d-1}$.  Put $t = \gamma(F_d) - k$, $1 \leq k \leq \gamma(F_d)$.  Then by Theorem \ref{th:mu=eps},
$$\mu_t(F_d) \leq \mu_t(Q_d) + 2^{d-1} =  \delta_t(Q_d) + 2^{d-1} = 4k + 2^{d-1}.$$
If $\mu_t(F_d) < 4k + 2^{d-1}$, then $\delta_t(F_d) < \delta_t(Q_d) + 2^{d-1}$ but the difference between the excesses is exactly $2^{d-1}$ by the above observation.
\end{proof}

\section{Discussion}

In summary, it appears that skewness is the more natural invariant and
we conjecture that 
the maximum difference between skewness $\mu$ and its (integer) Euler lower bound $\varepsilon$, over complete and bipartite complete graphs for all orientable surfaces, is {\it finite}.  
For cubes and folded cubes, it is zero.  

Circulant graphs are a natural common generalization of complete and regular bipartite complete graphs.  The planar excess is 4 for any triangle free, 4-regular graph. For circulant graphs $c(n,k) := C(n, \{1,k\})$ with $n$ vertices and jump-lengths $1, k < n/2$, we have $\delta_0(c(n,k)) = 4 = \mu_0(c(n,k))$ when $n = 2sk$ for $s, k \geq 2$ integers \cite{jko2025}.  
But Ho \cite{ho} proved $\nu_0(c(3k{+}1,k) {=} k{+}1$ for $k \geq 3$.  His lower-bound argument depends on the particular metric structure of this graph and uses induction on $k$ for a proof by contradiction.  Thus, for $c(3k{+}1,k)$, {\it the difference between planar excess and planar crossing number can be arbitrarily large}. Is this also true for the difference with {\it skewness}?  

Other pairs of invariants behave like $(\varepsilon, \mu)$ for suitable graph families.  For instance, (maximum degree $\Delta$, matching book thickness) \cite{alam-21, kjo-circ, ojk, ysl}, ($\lceil \Delta/2 \rceil$, linear arboricity) \cite{aeh}, (book thickness, star arboricity book thickness) \cite{pck-at}, inter alia.
These are similar to Vizing's theorem that the chromatic index cannot exceed the maximum degree by more than 1.  What is the extent of this ``Vizing Phenomenon'' and what is its theoretical significance?



\end{document}